\documentclass{article}
\usepackage{amssymb}
\usepackage{amsmath}
\usepackage{amsfonts}
\usepackage{hyperref}

\newtheorem{theorem}{Theorem}

\begin{document}

\title{Smoothness of radial solutions to Monge-Amp\`{e}re equations}
\author{Cristian Rios \\
%EndAName
University of Calgary\\
Calgary, Alberta \and Eric T. Sawyer \\
%EndAName
McMaster University\\
Hamilton, Ontario}
\maketitle

\section{Introduction}

It is well known that the radial homogeneous functions $u=c_{m,n}\left\vert
x\right\vert ^{2+\frac{2m}{n}}$ provide \emph{nonsmooth} solutions to the
Monge-Amp\`{e}re equation $\det D^{2}u=\left\vert x\right\vert ^{2m}$ with
smooth right hand side when $m\in \mathbb{N}\setminus n\mathbb{N}$. This
raises the question of when radial solutions $u$ to the generalized equation 
\begin{equation}
\det D^{2}u=k\left( x,u,Du\right) ,\ \ \ \ \ x\in \mathbb{B}_{n},
\label{DirMA}
\end{equation}%
are smooth, given that $k$ is smooth and nonnegative. When $u$ is radial, (%
\ref{DirMA}) reduces to a nonlinear ODE on $\left[ 0,1\right) $ that is
singular at the endpoint $0$. It is thus easy to prove that $u$ is always
smooth away from the origin, even where $k$ vanishes, but smoothness at the
origin is more complicated, and determined by the order of vanishing of $k$
there.

In fact, Monn \cite{Mo} proves that if $k=k\left( x\right) $ is independent
of $u$ and $Du$, then a radial solution $u$ to (\ref{DirMA}) is smooth if $%
k^{\frac{1}{n}}$ is smooth, and Derridj \cite{De} has extended this
criterion to the case when $k\left( x,u,Du\right) =f\left( \frac{\left\vert
x\right\vert ^{2}}{2},u,\frac{\left\vert \nabla u\right\vert ^{2}}{2}\right) 
$ factors as%
\begin{equation}
f\left( t,\xi ,\zeta \right) =\kappa \left( t\right) \phi \left( t,\xi
,\zeta \right)  \label{factors}
\end{equation}%
with $\kappa $ smooth and nonnegative on $\left[ 0,1\right) $, $\kappa
\left( 0\right) =0$, and $\phi $ smooth and positive on $\left[ 0,1\right)
\times \mathbb{R}\times \left[ 0,\infty \right) $. Moreover, Monn also shows
that $u$ is smooth if $k=k\left( x\right) $ vanishes to $\emph{infinite}$
order at the origin.

These results leave open the case when $k$\ has the general form $k\left(
x,u,Du\right) $ and vanishes to infinite order at the origin. The purpose of
this paper is to show that radial solutions $u$ are smooth in this remaining
case as well. The following theorem encompasses all of the afore-mentioned
results, and applies to generalized convex solutions $u$ and also with $%
f=\kappa \phi $ as in (\ref{factors}) but where $\phi $ is only assumed
positive and bounded, not smooth.

\begin{theorem}
\label{radialsmooth}Suppose that $u$ is a generalized convex radial solution
(in the sense of Alexandrov) to the generalized Monge-Amp\`{e}re equation (%
\ref{DirMA}) with%
\begin{equation*}
k\left( x,u,Du\right) =f\left( \frac{\left\vert x\right\vert ^{2}}{2},u,%
\frac{\left\vert \nabla u\right\vert ^{2}}{2}\right)
\end{equation*}%
where $f$ is smooth and nonnegative on $\left[ 0,1\right) \times \mathbb{R}%
\times \left[ 0,\infty \right) $. Then $u$ is smooth in the deleted ball $%
\mathbb{B}_{n}\setminus \left\{ 0\right\} $.\medskip \newline
Suppose moreover that there are positive constants $c,C$ such that%
\begin{equation}
cf\left( t,0,0\right) \leq f\left( t,\xi ,\zeta \right) \leq Cf\left(
t,0,0\right)  \label{factors'}
\end{equation}%
for $\left( \xi ,\zeta \right) $ near $\left( 0,0\right) $. Let $\tau \in 
\mathbb{Z}_{+}\cup \left\{ \infty \right\} $ be the order of vanishing of $%
f\left( t,0,0\right) $ at $0$. Then $u$ is smooth at the origin if and only
if $\tau \in n\mathbb{Z}_{+}\cup \left\{ \infty \right\} $.
\end{theorem}

The case when $k=k\left( x\right) $ is independent of $u$ and $Du$ is
handled by Monn in \cite{Mo} using an explicit formula for $u$ in terms of $%
k $:%
\begin{equation}
g\left( t\right) =C+\left( \frac{n}{2}\right) ^{\frac{1}{n}}\int_{0}^{t}%
\frac{\left( \int_{0}^{s}w^{\frac{n}{2}}f\left( w\right) \frac{dw}{w}\right)
^{\frac{1}{n}}}{\sqrt{s}}ds.  \label{radialform}
\end{equation}
where $u\left( x\right) =g\left( \frac{r^{2}}{2}\right) $ and $k\left(
x\right) =f\left( \frac{r^{2}}{2}\right) \geq 0$ with $r=\left\vert
x\right\vert $, $x\in \mathbb{R}^{n}$. In the case $k$ vanishes to infinite
order at the origin, an inequality of Hadamard is used as well. The
following scale invariant version follows from Corollary 5.2 in \cite{Mo}:%
\begin{equation}
\max_{0\leq t\leq x}\left\vert F^{\left( \ell \right) }\left( t\right)
\right\vert \leq C_{k,\ell }F\left( x\right) ^{\frac{k-\ell }{k}}\max_{0\leq
t\leq x}\left\vert F^{\left( k\right) }\left( t\right) \right\vert ^{\frac{%
\ell }{k}},\ \ \ \ \ 0\leq x\leq 1,  \label{mon}
\end{equation}%
for all $1\leq \ell \leq k-1$ and $k\in \mathbb{N}$ provided $F$ is smooth,
nondecreasing on $\left[ 0,1\right) $ and vanishes to infinite order at $0$.

\section{Proof of Theorem \protect\ref{radialsmooth}}

We begin by considering Theorem \ref{radialsmooth} in the case that $u$ is a
classical $C^{2}$ solution to (\ref{DirMA}) and $f$ satisfies (\ref{factors}%
) where $f\left( t,0,0\right) $ vanishes to \emph{finite} order $\ell $ at $%
0 $. If $k$ is independent of $u$ and $Du$, Monn uses formula (\ref%
{radialform}) in \cite{Mo} to show that $u$ is smooth when $f\left( w\right)
^{\frac{1}{n}}$ is smooth. In particular this applies when $\ell \in n%
\mathbb{Z}_{+}$. In the general case, we note that (\ref{factors'}) implies (%
\ref{factors}), the assumption made in \cite{De}. Indeed, using $f^{\left(
k\right) }\left( 0,\xi ,\zeta \right) =0$ for $0\leq k\leq \ell -1$ we can
write%
\begin{equation*}
f\left( s,\xi ,\zeta \right) =\int_{0}^{1}\frac{\left( 1-t\right) ^{\ell -1}%
}{\left( \ell -1\right) !}\frac{d^{\ell }}{dt^{\ell }}f\left( ts,\xi ,\zeta
\right) dt=s^{\ell }\psi \left( s,\xi ,\zeta \right) ,
\end{equation*}%
where $\psi \left( s,\xi ,\zeta \right) $ is smooth and $\psi \left( 0,\xi
,\zeta \right) =\frac{f^{\left( \ell \right) }\left( 0,\xi ,\zeta \right) }{%
\ell !}>0$. Thus the results of Derridj \cite{De} apply to show that $u$ is
smooth for general $k$ when $\ell \in n\mathbb{Z}_{+}$.

\subsection{Generalized Monge-Amp\`{e}re equations}

We now consider radial \emph{generalized convex} solutions $u$ to the
generalized Monge-Amp\`{e}re equation (\ref{DirMA}) where we assume $k\left(
\cdot ,u,q\right) $ and $k\left( x,u,\cdot \right) $ are radial. We first
establish that $u\in C^{2}\left( \mathbb{B}_{n}\right) \cap C^{\infty
}\left( \mathbb{B}_{n}\setminus \left\{ 0\right\} \right) $. We note that
results of Guan, Trudinger and Wang in \cite{Gu2} and \cite{GuTrWa} yield $%
u\in C^{1,1}\left( \mathbb{B}_{n}\right) $ for many $k$ in (\ref{DirMA}),
but not in the generality possible in the radial case here. In order to deal
with general $k$ it would be helpful to have a formula for $u$ in terms of $%
k $, but this is problematic. Instead we prove Theorem \ref{radialsmooth}
for general $k$ \emph{without} solving for the solution explicitly, but
using an inductive argument that is based on Lemma \ref{monot} when $k$
vanishes to infinite order at the origin.

Assume that $u$ is a generalized convex solution of (\ref{DirMA}) in the
sense of Alexandrov (see \cite{Al} and \cite{ChYa}) and define $\varphi
\left( t\right) $ by%
\begin{equation}
\varphi \left( \frac{r^{2}}{2}\right) =k\left( x,u\left( x\right) ,Du\left(
x\right) \right) =f\left( \frac{\left\vert x\right\vert ^{2}}{2},u\left(
x\right) ,\frac{\left\vert \nabla u\left( x\right) \right\vert ^{2}}{2}%
\right) .  \label{name}
\end{equation}%
Then $\varphi $ is bounded since $u$ is Lipschitz continuous. It follows
that the \emph{convex radial} function $u$ is continuously differentiable at
the origin, since otherwise it would have a conical singularity there and
its representing measure $\mu _{u}$ would have a Dirac component at the
origin. Let $g$ be given by formula (\ref{radialform}) with $\varphi $ in
place of $f$, i.e.%
\begin{equation}
g\left( t\right) =C_{u}+\left( \frac{n}{2}\right) ^{\frac{1}{n}}\int_{0}^{t}%
\frac{\left( \int_{0}^{s}w^{\frac{n}{2}}\varphi \left( w\right) \frac{dw}{w}%
\right) ^{\frac{1}{n}}}{\sqrt{s}}ds,  \label{gis}
\end{equation}
and with constant $C_{u}$ chosen so that $u$ and $\widetilde{u}$ agree on
the unit sphere where%
\begin{equation}
\widetilde{u}\left( x\right) =g\left( \frac{r^{2}}{2}\right) ,\ \ \ \ \
0\leq r<1.  \label{utilda}
\end{equation}

We claim that $\widetilde{u}$ is a generalized convex solution to (\ref%
{DirMA}) in the sense of Alexandrov. To see this we first note that $D^{2}%
\widetilde{u}\left( r\mathbf{e}_{1}\right) =\left[ 
\begin{array}{llll}
g^{\prime \prime }r^{2}+g^{\prime } & 0 & \cdots  & 0 \\ 
0 & g^{\prime } & \cdots  & 0 \\ 
\vdots  & \vdots  & \ddots  & \vdots  \\ 
0 & 0 & \cdots  & g^{\prime }%
\end{array}%
\right] $ is positive semidefinite, hence $\widetilde{u}$ is convex. To
prove that the representing measure $\mu _{\widetilde{u}}$ of $\widetilde{u}$%
\ is $kdx$ it suffices to show, since both $g$ and $f$ are radial, that%
\begin{equation*}
\mu _{\widetilde{u}}\left( E\right) =\left\vert B_{\widetilde{u}}\left(
E\right) \right\vert =\int_{E}k
\end{equation*}%
for all annuli $E=\left\{ x\in \mathbb{B}_{n}:r_{1}<\left\vert x\right\vert
<r_{2}\right\} $, $0<r_{1}<r_{2}<1$ where%
\begin{equation*}
B_{\widetilde{u}}\left( E\right) =\cup _{r_{1}<\left\vert x\right\vert
<r_{2}}\left\{ \nabla \widetilde{u}_{1}\left( x\right) \right\} =\left\{
a\in \mathbb{B}_{n}:\frac{\partial }{\partial r}\widetilde{u}\left( r_{1}%
\mathbf{e}_{1}\right) <\left\vert a\right\vert <\frac{\partial }{\partial r}%
\widetilde{u}\left( r_{2}\mathbf{e}_{1}\right) \right\} .
\end{equation*}%
Since $\frac{\partial }{\partial r}\widetilde{u}\left( r_{i}\mathbf{e}%
_{i}\right) =g^{\prime }\left( \frac{r_{i}^{2}}{2}\right) r_{i}=g^{\prime
}\left( t_{i}\right) \sqrt{2t_{i}}$ with $t_{i}=\frac{r_{i}^{2}}{2}$, we
thus have%
\begin{eqnarray*}
\left\vert B_{\widetilde{u}}\left( E\right) \right\vert  &=&\left\vert
\left\{ a\in \mathbb{B}_{n}:g^{\prime }\left( t_{1}\right) \sqrt{2t_{1}}%
<\left\vert a\right\vert <g^{\prime }\left( t_{2}\right) \sqrt{2t_{2}}%
\right\} \right\vert  \\
&=&\frac{\omega _{n}}{n}\left\{ g^{\prime }\left( t_{2}\right) ^{n}\left(
2t_{2}\right) ^{\frac{n}{2}}-g^{\prime }\left( t_{1}\right) ^{n}\left(
2t_{1}\right) ^{\frac{n}{2}}\right\}  \\
&=&\frac{\omega _{n}}{n}\frac{n}{2}2^{\frac{n}{2}}\int_{t_{1}}^{t_{2}}w^{%
\frac{n}{2}}f\left( w\right) \frac{dw}{w} \\
&=&\omega _{n}\int_{r_{1}}^{r_{2}}r^{n-1}\varphi \left( \frac{r^{2}}{2}%
\right) dr=\int_{E}k.
\end{eqnarray*}%
In particular the \emph{convex radial} function $\widetilde{u}$ must be
continuously differentiable, since otherwise there is a jump discontinuity
in the radial derivative of $\widetilde{u}$ at some distance $r$ from the
origin that results in a singular component in $\mu _{\widetilde{u}}$
supported on the sphere of radius $r$.

Now uniqueness of Alexandrov solutions to the Dirichlet problem (see e.g. 
\cite{ChYa}) yields $u=\widetilde{u}$, and hence $u\in C^{1}\left( \mathbb{B}%
_{n}\right) $. Thus $\varphi \in C\left[ 0,1\right) $ and from (\ref{utilda}%
) we have $u\left( x\right) =g\left( \frac{\left\vert x\right\vert ^{2}}{2}%
\right) $ and%
\begin{equation}
\varphi \left( t\right) =f\left( t,g\left( t\right) ,tg^{\prime }\left(
t\right) ^{2}\right) ,  \label{phiis}
\end{equation}%
where using (\ref{gis}) we compute that%
\begin{equation}
g^{\prime }\left( t\right) =\left\{ \frac{n}{2}t^{-\frac{n}{2}%
}\int_{0}^{t}s^{\frac{n}{2}-1}\varphi \left( s\right) ds\right\} ^{\frac{1}{n%
}}.  \label{g'}
\end{equation}%
In particular $g^{\prime }\in C\left[ 0,1\right) $. We now obtain by
induction that $g\in C^{\infty }\left( 0,1\right) $, hence $u\in C^{\infty
}\left( \mathbb{B}_{n}\setminus \left\{ 0\right\} \right) $. Indeed, if $%
g\in C^{\ell }\left( 0,1\right) $\ then (\ref{phiis}) implies $\varphi \in
C^{\ell -1}\left( 0,1\right) $ and then (\ref{gis}) implies $g\in C^{\ell
+1}\left( 0,1\right) $.

It will be convenient to use fractional integral operators at this point.
For $\beta >0$ and $f$ continuous define%
\begin{eqnarray*}
T_{\beta }f\left( s\right) &=&\int_{0}^{s}\left( \frac{w}{s}\right) ^{\beta
}f\left( w\right) \frac{dw}{w},\ \ \ \ \ s\neq 0, \\
T_{\beta }f\left( 0\right) &=&\frac{1}{\beta }f\left( 0\right) ,
\end{eqnarray*}%
so that%
\begin{equation}
g\left( t\right) =C+\left( \frac{n}{2}\right) ^{\frac{1}{n}%
}\int_{0}^{t}\left( T_{\frac{n}{2}}f\left( s\right) \right) ^{\frac{1}{n}}ds.
\label{gequals}
\end{equation}%
We claim that for $f$ smooth, nonnegative and of finite type $\ell $, $\ell
\in \mathbb{Z}_{+}$, the same is true of $T_{\beta }f$ for all $\beta >0$.
This follows immediately from the identity%
\begin{equation}
\frac{d^{k}}{ds^{k}}T_{\beta }f\left( s\right) =T_{\beta +k}f^{\left(
k\right) }\left( s\right) ,\ \ \ \ \ k\in \mathbb{N},  \label{iden}
\end{equation}%
and the estimate%
\begin{equation*}
T_{\beta +k}f^{\left( k\right) }\left( s\right) =\frac{1}{\beta +k}f^{\left(
k\right) }\left( 0\right) +O\left( \left\vert s\right\vert \right) .
\end{equation*}%
When $k=1$, (\ref{iden}) follows from differentiating and then integrating
by parts, and the general case is then obtained by iteration.

Now suppose that $f$ satisfies (\ref{factors'}) and let%
\begin{equation*}
\kappa \left( t\right) =f\left( t,0,0\right)
\end{equation*}%
vanish to infinite order at $0$. If $\kappa $ vanishes in a neighbourhood of 
$0$ then so does $g$ and we have $g\in C^{\infty }\left[ 0,1\right) $ and $%
u\in C^{\infty }\left( \mathbb{B}_{n}\right) $. Thus we will assume $%
\int_{0}^{t}\kappa >0$ for $t>0$ in what follows. Note that (\ref{iden})
then implies that $T_{\frac{n}{2}}\kappa \left( t\right) $ is smooth and
positive on $\left( 0,1\right) $ and vanishes to infinite order at $0$.
Since $g^{\prime }\in C\left[ 0,1\right) $, it follows that $\varphi \left(
t\right) \leq C\kappa \left( t\right) $. Thus we have the inequality $T_{%
\frac{n}{2}}\varphi \left( t\right) \leq CT_{\frac{n}{2}}\kappa \left(
t\right) $, and from (\ref{g'}) we now conclude that $g^{\prime }\left(
t\right) $ also vanishes to infinite order at $0$. Now $\varphi \left(
t\right) \approx \kappa \left( t\right) $ from (\ref{factors'}), and so also 
$T_{\frac{n}{2}}\varphi \left( t\right) \approx T_{\frac{n}{2}}\kappa \left(
t\right) $. From%
\begin{equation}
g^{\prime \prime }\left( t\right) =\frac{\varphi \left( t\right) }{2t\left( 
\frac{n}{2}T_{\frac{n}{2}}\varphi \left( t\right) \right) ^{1-\frac{1}{n}}}-%
\frac{1}{2t}\left( \frac{n}{2}T_{\frac{n}{2}}\varphi \left( t\right) \right)
^{\frac{1}{n}},  \label{g''}
\end{equation}%
we then have%
\begin{equation}
\left\vert g^{\prime \prime }\left( t\right) \right\vert \leq C\frac{\kappa
\left( t\right) }{2t\left( \frac{n}{2}T_{\frac{n}{2}}\kappa \left( t\right)
\right) ^{1-\frac{1}{n}}}+C\frac{1}{2t}\left( \frac{n}{2}T_{\frac{n}{2}%
}\kappa \left( t\right) \right) ^{\frac{1}{n}},\ \ \ \ \ 0<t<1.
\label{approx}
\end{equation}

An application of (\ref{mon}) with $\ell =1$, $k>n$ and $F\left( t\right)
=\int_{0}^{t}s^{\frac{n}{2}-1}\kappa \left( s\right) ds$ yields $t^{\frac{n}{%
2}}\kappa \left( t\right) =F^{\prime }\left( t\right) \leq CF\left( t\right)
^{1-\frac{1}{k}}$ and so the first term on the right side of (\ref{approx})
is bounded by a multiple of $t^{-\frac{1}{2}}F\left( t\right) ^{\frac{1}{n}-%
\frac{1}{k}}$. Thus the right side of (\ref{approx}), and hence also $%
g^{\prime \prime }\left( t\right) $, vanishes to infinite order at $0$. In
particular $g^{\prime \prime }\in C\left[ 0,1\right) $ and we conclude $u\in
C^{2}\left( \mathbb{B}_{n}\right) $ in this case as well.

Summarizing, we have $u\in C^{\infty }\left( \mathbb{B}_{n}\setminus \left\{
0\right\} \right) $, and in the case $f$ satifies (\ref{factors'}), we also
have $u\in C^{2}\left( \mathbb{B}_{n}\right) $. Thus from above we have that 
\begin{equation*}
\varphi \left( t\right) =f\left( t,g\left( t\right) ,tg^{\prime }\left(
t\right) ^{2}\right) =\kappa \left( t\right) \phi \left( t,g\left( t\right)
,tg^{\prime }\left( t\right) ^{2}\right) ,
\end{equation*}%
where $u\left( x\right) =g\left( \frac{\left\vert x\right\vert ^{2}}{2}%
\right) \in C^{2}\left( \mathbb{B}_{n}\right) $, $g$ is given by (\ref{gis})
and $\varphi \in C^{1}\left[ 0,1\right) $ by (\ref{name}). Note that we
cannot use (\ref{mon}) on the function $\int_{0}^{t}s^{\frac{n}{2}-1}\varphi
\left( s\right) ds$ here since we have no \textit{a priori} control on
higher derivatives of $\varphi \left( s\right) =f\left( s,g\left( s\right)
,sg^{\prime }\left( s\right) ^{2}\right) $. Instead we will use (\ref{mon})
on the function $\int_{0}^{t}s^{\frac{n}{2}-1}\kappa \left( s\right) ds$
together with an inductive argument to control derivatives of $g$.

From above we have that $g^{\prime \prime }\in C\left[ 0,1\right) \cap
C^{\infty }\left( 0,1\right) $. Now differentiate (\ref{g''}) for $t>0$
using (\ref{iden}) to obtain%
\begin{eqnarray}
&&g^{\prime \prime \prime }\left( t\right)  \label{g'''} \\
&&\ \ \ \ \ =\frac{1}{2}\left( \frac{n}{2}\right) ^{\frac{1}{n}-1}\left\{ 
\frac{\varphi ^{\prime }\left( t\right) }{tT_{\frac{n}{2}}\varphi \left(
t\right) ^{1-\frac{1}{n}}}-\left( \frac{1}{n}-1\right) \frac{\varphi \left(
t\right) T_{\frac{n}{2}+1}\varphi ^{\prime }\left( t\right) }{tT_{\frac{n}{2}%
}\varphi \left( t\right) ^{2-\frac{1}{n}}}-\frac{\varphi \left( t\right) }{%
t^{2}T_{\frac{n}{2}}\varphi \left( t\right) ^{1-\frac{1}{n}}}\right\}  \notag
\\
&&\ \ \ \ \ -\frac{1}{2}\left( \frac{n}{2}\right) ^{\frac{1}{n}}\left\{ 
\frac{1}{n}\frac{T_{\frac{n}{2}+1}\varphi ^{\prime }\left( t\right) }{tT_{%
\frac{n}{2}}\varphi \left( t\right) ^{1-\frac{1}{n}}}-\frac{T_{\frac{n}{2}%
}\varphi \left( t\right) ^{\frac{1}{n}}}{t^{2}}\right\} ,  \notag
\end{eqnarray}%
and then compute that%
\begin{eqnarray}
\varphi ^{\prime }\left( t\right) &=&\kappa ^{\prime }\left( t\right) \phi
\left( t,g\left( t\right) ,tg^{\prime }\left( t\right) ^{2}\right)
\label{phi'} \\
&&+\kappa \left( t\right) \phi _{1}\left( t,g\left( t\right) ,tg^{\prime
}\left( t\right) ^{2}\right)  \notag \\
&&+\kappa \left( t\right) \phi _{2}\left( t,g\left( t\right) ,tg^{\prime
}\left( t\right) ^{2}\right) g^{\prime }\left( t\right)  \notag \\
&&+\kappa \left( t\right) \phi _{3}\left( t,g\left( t\right) ,tg^{\prime
}\left( t\right) ^{2}\right) \left\{ g^{\prime }\left( t\right)
^{2}+2tg^{\prime }\left( t\right) g^{\prime \prime }\left( t\right) \right\}
.  \notag
\end{eqnarray}%
We will now use $\varphi \approx \kappa $, (\ref{g'''}), (\ref{phi'}) and (%
\ref{mon}) applied with $F\left( t\right) =\int_{0}^{t}s^{\frac{n}{2}%
-1}\kappa \left( s\right) ds$, to show that $g^{\prime \prime \prime }$
vanishes to infinite order at $0$ and $g^{\prime \prime \prime }\in C\left[
0,1\right) $.

To see this, we first note that $F$ is smooth, nonnegative and vanishes to
infinite order at $0$ since the same is true of $\kappa $. Next, for any $%
\ell \geq 1$ and $\varepsilon >0$,\ (\ref{mon}) with $k$ large enough yields%
\begin{equation}
\sup_{0<s\leq t}\left\vert F^{\left( \ell \right) }\left( s\right)
\right\vert \leq C_{\varepsilon ,\ell }F\left( t\right) ^{1-\varepsilon }.
\label{mon'}
\end{equation}%
Moreover we have%
\begin{align}
\left\vert \beta T_{\beta }h\left( t\right) \right\vert & \leq \sup_{0<s\leq
t}\left\vert h\left( s\right) \right\vert ,  \label{more} \\
F\left( t\right) & =t^{\frac{n}{2}}T_{\frac{n}{2}}\kappa \left( t\right) , 
\notag \\
T_{\frac{n}{2}}\varphi \left( t\right) & \approx T_{\frac{n}{2}}\kappa
\left( t\right) .  \notag
\end{align}%
Now using 
\begin{eqnarray*}
F^{\prime }\left( t\right) &=&t^{\frac{n}{2}-1}\kappa \left( t\right) , \\
F^{\prime \prime }\left( t\right) &=&t^{\frac{n}{2}-1}\kappa ^{\prime
}\left( t\right) +\left( \frac{n}{2}-1\right) t^{\frac{n}{2}-2}\kappa \left(
t\right) ,
\end{eqnarray*}%
yields%
\begin{equation*}
\left\vert \kappa ^{\prime }\left( t\right) \phi \left( t,g\left( t\right)
,tg^{\prime }\left( t\right) ^{2}\right) \right\vert \leq C\left\vert \kappa
^{\prime }\left( t\right) \right\vert =C\left\vert t^{1-\frac{n}{2}%
}F^{\prime \prime }\left( t\right) -\left( \frac{n}{2}-1\right) t^{-\frac{n}{%
2}}F^{\prime }\left( t\right) \right\vert ,
\end{equation*}%
and an application of (\ref{mon'}) gives 
\begin{equation*}
\left\vert \kappa ^{\prime }\left( t\right) \phi \left( t,g\left( t\right)
,tg^{\prime }\left( t\right) ^{2}\right) \right\vert \leq C_{\varepsilon
}t^{-\frac{n}{2}}F\left( t\right) ^{1-\varepsilon }.
\end{equation*}%
We obtain similar estimates for the remaining terms in (\ref{phi'}) and
altogether this yields%
\begin{equation*}
\left\vert \varphi ^{\prime }\left( t\right) \right\vert \leq C_{\varepsilon
}t^{-\alpha }F\left( t\right) ^{1-\varepsilon },\ \ \ \ \ \text{for some }%
\alpha >0.
\end{equation*}%
Using the second and third lines in (\ref{more}) now shows that the first
term in braces in (\ref{g'''}) satisfies%
\begin{equation*}
\left\vert \frac{\varphi ^{\prime }\left( t\right) }{tT_{\frac{n}{2}}\varphi
\left( t\right) ^{1-\frac{1}{n}}}\right\vert \leq C_{\varepsilon }\frac{%
t^{-\alpha }F\left( t\right) ^{1-\varepsilon }}{tT_{\frac{n}{2}}\varphi
\left( t\right) ^{1-\frac{1}{n}}}\approx C_{\varepsilon }t^{\frac{n}{2}%
\left( 1-\varepsilon \right) -\alpha -1}T_{\frac{n}{2}}\kappa \left(
t\right) ^{\frac{1}{n}-\varepsilon },
\end{equation*}%
which vanishes to infinite order at $0$ if $0<\varepsilon <\frac{1}{n}$.
Similar arguments, using (\ref{phi'}) and the first line in (\ref{more}) to
estimate $T_{\frac{n}{2}+1}\varphi ^{\prime }\left( t\right) $, apply to the
remaining terms in (\ref{g'''}), and this completes the proof that $%
g^{\prime \prime \prime }$ vanishes to infinite order at $0$ and $g^{\prime
\prime \prime }\in C\left[ 0,1\right) $.

We now observe that we can

\begin{itemize}
\item continue to differentiate (\ref{g'''}) to obtain a formula for $%
g^{\left( \ell \right) }$ involving only appropriate powers of $T_{\frac{n}{2%
}}\varphi \left( t\right) \approx T_{\frac{n}{2}}\kappa \left( t\right) $ in
the denominator, and derivatives of $\varphi $ of order at most $\ell -2$ in
the numerator,

\item and continue to differentiate (\ref{phi'}) to obtain a formula for $%
\varphi ^{\left( \ell -2\right) }$ involving derivatives of $g$ of order at
most $\ell -1$.
\end{itemize}

It is now clear that the above arguments apply to prove that derivatives of $%
g\left( t\right) $ of all orders vanish to infinite order at $0$ and are
continuous on $\left[ 0,1\right) $. This shows that $g$ is smooth on $\left[
0,1\right) $ and thus that $u$ is smooth on $\mathbb{B}_{n}$.

\end{document}